
\def\title#1{{\titlefont\leftskip-.5em\noindent #1\bigskip}}

\def\author#1{{\largefont\noindent #1}\medskip}

\def\beginlinemode{\endmode
 \begingroup\obeylines\def\endmode{\par\endgroup}}
\let\endmode=\par

\newbox\theaddress
\def\address{\smallskip\beginlinemode\parindent 0in\getaddress}
{\obeylines
\gdef\getaddress #1 
 #2
 {#1\gdef\addressee{#2}%
   \global\setbox\theaddress=\vbox\bgroup\raggedright%
    \everypar{\hangindent2em}#2
   \def\endaddress{\egroup\endgroup \copy\theaddress \medskip}}}

\def\thanks#1{\footnote{}{\eightpoint #1}}

\long\def\Abstract#1{{\eightpoint\narrower\vskip\baselineskip\noindent
#1\smallskip}}

\def\skipfirstword#1 {}

\def\ir#1{\csname #1\endcsname}

\newdimen\currentht
\newbox\droppedletter
\newdimen\droppedletterwdth
\newdimen\drophtinpts
\newdimen\dropindent

\def\irrnSection#1#2{\edef\tttempcs{\ir{#2}}
\currentht-\pagetotal\advance\currentht by-\ht\footins
\advance\currentht by\vsize
\ifdim\currentht<1.5in\par\vfill\eject\else\vbox to.25in{\vfil}\fi
{\largefont\noindent{}\hskip-.5em\expandafter\skipfirstword\tttempcs. #1}
\vskip\baselineskip\noindent\drop }

\def\boxit#1{\vbadness10000\vbox{\hrule\hbox{\vrule\kern3pt\vbox{\kern3pt
          #1\kern3pt}\kern3pt\vrule}\hrule}}

\def\drop#1{\setbox\droppedletter=\hbox{\boxit{\hbox{\dropfont #1}}}%
    \droppedletterwdth=\wd\droppedletter
    \advance\droppedletterwdth by.25em%
    \drophtinpts\normalbaselineskip
    \advance\drophtinpts by-\dp\droppedletter
    \advance\drophtinpts by.5ex
    \global\dropindent\droppedletterwdth
    \global\advance\dropindent by-.5em
    \global\hangindent\dropindent
    \global\hangafter-2
    \setbox\droppedletter=\hbox 
                to\droppedletterwdth{\raise-\drophtinpts\box\droppedletter\hfil}%
    \dp\droppedletter 0in%
    \ht\droppedletter 0in\llap{\box\droppedletter}}

\def\irSubsection#1#2{\edef\tttempcs{\ir{#2}}
\vskip\baselineskip\penalty-3000
{\bf\noindent \expandafter\skipfirstword\tttempcs. #1}
\vskip\baselineskip}

\def\irSubsubsection#1#2{\edef\tttempcs{\ir{#2}}
\vskip\baselineskip\penalty-3000
{\bf\noindent \expandafter\skipfirstword\tttempcs. #1}
\vskip\baselineskip}

\def\References{\vbox to.25in{\vfil}\noindent{}\hskip-.5em{\bf References}
\vskip\baselineskip\par}

\def\baselinebreak{\par \ifdim\lastskip<\baselineskip
         \removelastskip\penalty-200\vskip\baselineskip\fi}

\long\def\prclm#1#2#3{\baselinebreak
\noindent{\bf \csname #2\endcsname}:\enspace{\sl #3\par}\baselinebreak}

\def\Prf{\noindent{\bf Proof}: }

\def\rem#1#2{\baselinebreak\noindent{\bf \csname #2\endcsname}:\enspace }

\def\qed{{\hfill$\diamondsuit$}\vskip\baselineskip}

\def\bibitem#1{\par\indent\llap{\rlap{\bf [#1]}\indent}\indent\hangindent
2\parindent\ignorespaces}

\long\def\eatit#1{}

\def\leftheadlinetext{}
\def\rightheadlinetext{}

\def\leftheadline{{\eightrm\folio\hfil \leftheadlinetext\hfil}}
\def\rightheadline{{\eightrm\hfil\rightheadlinetext\hfil\folio}}

\headline={\ifnum\pageno=1\hfil\else
\ifodd\pageno\rightheadline\else\leftheadline\fi\fi}

\def\tenpoint{\def\rm{\fam0\tenrm}
\textfont0=\tenrm \scriptfont0=\sevenrm \scriptscriptfont0=\fiverm
\textfont1=\teni \scriptfont1=\seveni \scriptscriptfont1=\fivei
\def\mit{\fam1} \def\oldstyle{\fam1\teni}
\textfont2=\tensy \scriptfont2=\sevensy \scriptscriptfont2=\fivesy
\def\cal{\fam2}
\textfont3=\tenex \scriptfont3=\tenex \scriptscriptfont3=\tenex
\def\it{\fam\itfam\tenit} 
\textfont\itfam=\tenit
\def\sl{\fam\slfam\tensl} 
\textfont\slfam=\tensl
\def\bf{\fam\bffam\tenbf} 
\textfont\bffam=\tenbf \scriptfont\bffam=\sevenbf
\scriptscriptfont\bffam=\fivebf
\def\tt{\fam\ttfam\tentt} 
\textfont\ttfam=\tentt
\normalbaselineskip=12pt
\setbox\strutbox=\hbox{\vrule height8.5pt depth3.5pt  width0pt}%
\normalbaselines\rm}

\def\eightpoint{\def\rm{\fam0\eightrm}%
\textfont0=\eightrm \scriptfont0=\sixrm \scriptscriptfont0=\fiverm
\textfont1=\eighti \scriptfont1=\sixi \scriptscriptfont1=\fivei
\def\mit{\fam1} \def\oldstyle{\fam1\eighti}%
\textfont2=\eightsy \scriptfont2=\sixsy \scriptscriptfont2=\fivesy
\def\cal{\fam2}%
\textfont3=\tenex \scriptfont3=\tenex \scriptscriptfont3=\tenex
\def\it{\fam\itfam\eightit} 
\textfont\itfam=\eightit
\def\sl{\fam\slfam\eightsl} 
\textfont\slfam=\eightsl
\def\bf{\fam\bffam\eightbf} 
\textfont\bffam=\eightbf \scriptfont\bffam=\sixbf
\scriptscriptfont\bffam=\fivebf
\def\tt{\fam\ttfam\eighttt} 
\textfont\ttfam=\eighttt
\normalbaselineskip=9pt%
\setbox\strutbox=\hbox{\vrule height7pt depth2pt  width0pt}%
\normalbaselines\rm}

\def\largefont{\def\rm{\fam0\largerm}
\textfont0=\largerm \scriptfont0=\largescriptrm \scriptscriptfont0=\tenrm
\textfont1=\largei \scriptfont1=\largescripti \scriptscriptfont1=\teni
\def\mit{\fam1} \def\oldstyle{\fam1\teni}
\textfont2=\largesy 
\def\cal{\fam2}
\def\it{\fam\itfam\largeit} 
\textfont\itfam=\largeit
\def\sl{\fam\slfam\largesl} 
\textfont\slfam=\largesl
\def\bf{\fam\bffam\largebf} 
\textfont\bffam=\largebf 
\def\tt{\fam\ttfam\largett} 
\textfont\ttfam=\largett
\normalbaselineskip=17.28pt
\setbox\strutbox=\hbox{\vrule height12.25pt depth5pt  width0pt}%
\normalbaselines\rm}

\def\titlefont{\def\rm{\fam0\titlerm}
\textfont0=\titlerm \scriptfont0=\largescriptrm \scriptscriptfont0=\tenrm
\textfont1=\titlei \scriptfont1=\largescripti \scriptscriptfont1=\teni
\def\mit{\fam1} \def\oldstyle{\fam1\teni}
\textfont2=\titlesy 
\def\cal{\fam2}
\def\it{\fam\itfam\titleit} 
\textfont\itfam=\titleit
\def\sl{\fam\slfam\titlesl} 
\textfont\slfam=\titlesl
\def\bf{\fam\bffam\titlebf} 
\textfont\bffam=\titlebf 
\def\tt{\fam\ttfam\titlett} 
\textfont\ttfam=\titlett
\normalbaselineskip=24.8832pt
\setbox\strutbox=\hbox{\vrule height12.25pt depth5pt  width0pt}%
\normalbaselines\rm}

\nopagenumbers

\font\eightrm=cmr8
\font\eighti=cmmi8
\font\eightsy=cmsy8
\font\eightbf=cmbx8
\font\eighttt=cmtt8
\font\eightit=cmti8
\font\eightsl=cmsl8
\font\sixrm=cmr6
\font\sixi=cmmi6
\font\sixsy=cmsy6
\font\sixbf=cmbx6

\font\largerm=cmr12 at 17.28pt
\font\largei=cmmi12 at 17.28pt
\font\largescriptrm=cmr12 at 14.4pt
\font\largescripti=cmmi12 at 14.4pt
\font\largesy=cmsy10 at 17.28pt
\font\largebf=cmbx12 at 17.28pt
\font\largett=cmtt12 at 17.28pt
\font\largeit=cmti12 at 17.28pt
\font\largesl=cmsl12 at 17.28pt

\font\titlerm=cmr12 at 24.8832pt
\font\titlei=cmmi12 at 24.8832pt
\font\titlesy=cmsy10 at 24.8832pt
\font\titlebf=cmbx12 at 24.8832pt
\font\titlett=cmtt12 at 24.8832pt
\font\titleit=cmti12 at 24.8832pt
\font\titlesl=cmsl12 at 24.8832pt
\font\dropfont=cmr12 at 24.8832pt

\tenpoint



\hsize 6.5in
\vsize 9.2in

\tolerance 3000
\hbadness 3000

\def\trans{H1}
\def\vanc{H2}
\def\ravello{H3}
\def\antican{H4}
\def\fatpts{H5}

\def\binom#1#2{\hbox{$\left(\matrix{#1\cr #2\cr}\right)$}}
\def\C#1{\hbox{$\cal #1$}}
\def\s{\hbox{$\scriptstyle\cal S$}}
\def\r{\hbox{$\scriptstyle\cal R$}}

\def\pr#1{\hbox{{\bf P}${}^{#1}$}}
\def\leftheadlinetext{Brian Harbourne}
\def\rightheadlinetext{Generators for Symbolic Powers}

\title{Generators for Symbolic Powers of 
Ideals Defining General Points of \pr2}

\author{Brian Harbourne}

\address
Department of Mathematics and Statistics
University of Nebraska-Lincoln
Lincoln, NE 68588-0323
email: bharbourne@unl.edu
WEB: http://www.math.unl.edu/$\sim$bharbour/
\smallskip
August 6, 1995\endaddress
\vskip-\baselineskip

\thanks{\vskip -6pt
\noindent This work benefitted from a National Science Foundation grant.
\smallskip
\noindent 1980 {\it Mathematics Subject Classification. } Primary 13P10, 14C99. 
Secondary 13D02, 13H15.
\smallskip
\noindent {\it Key words and phrases. }  Ideal generation conjecture, symbolic powers,
fat points.\smallskip}

\vskip-6pt
\Abstract{}

\irrnSection{Introduction}{intro}
Let $I$ be an ideal,
homogeneous with respect to the usual grading,
in a polynomial ring
$R=k[x_0,\ldots,x_n]$ in $n+1$ variables (over an algebraically closed 
field $k$). Denote the graded component of $I$
of degree $d$ by $I_d$, and likewise the $k$-vector space of homogeneous
forms of $R$ of degree $d$ by $R_d$. Since $I$ is a graded $R$-module,
we have $k$-linear maps $\mu_{d,i}:I_d\otimes R_i\to I_{d+i}$ given
for each $i$ and $d$ by multiplication; when the index $i$ (or the indices $i$ and $d$)
are clear from context, we will denote $\mu_{d,i}$ by $\mu_d$ or simply $\mu$.

\rem{Definition}{mrp} Say that $I$ has the {\it maximal rank property\/}
if for each $d$ either the kernel or cokernel
of the homomorphism $\mu:I_d\otimes R_1\to I_{d+1}$ vanishes. \qed

Note that the number $\nu_{d+1}$ of generators of $I$ of 
degree $d+1$ in a minimal homogeneous
set of generators is just the $k$-dimension of the cokernel of $\mu_{d,1}$.
Thus for an ideal $I$ with the maximal rank property
one only needs to know the Hilbert function $H_I(d)=\hbox{dim}(I_d)$
of $I$ to determine $\nu_d$ for every $d$. This suggests
why one might want to know under what conditions an ideal
has the maximal rank property.

One situation which has attracted interest is that of
points in projective space. If $p_1,\ldots,p_r\in\pr n$ are
distinct points, let $Z=p_1+\cdots+p_r$ denote the smooth subscheme
given by the union of the points. Its corresponding homogeneous ideal $I(Z)\subset R$
is $I(Z)=I(p_1)\cap\cdots\cap I(p_r)$, where $I(p_j)$ is for each $j$ the ideal 
generated by all forms vanishing at $p_j$. Concerning such ideals there is
the Ideal Generation Conjecture (IGC) of [Ro], [GO] and [GGR]:

\prclm{Ideal Generation Conjecture}{IGC}{The ideal $I(Z)$ has the maximal rank property 
for a general set $Z=p_1+\cdots+p_r$ of $r$ points in \pr n.}

Since, in the situation of the conjecture, the Hilbert function $H_{I(Z)}$ 
is known, the conjecture would allow one to determine numbers $\nu_d$ of generators
in each degree $d$, and indeed the conjecture
has been verified for a number of values of $r$ and $n$ (see [B], [GM], [HSV],
[L], [O], [Ra]). Much less is understood
or even conjectured in the more general situation of fat points
(although we would draw your attention to [Cat] and [I]):
a {\it fat point\/} subscheme $Z\subset \pr n$ is a subscheme
defined by a homogeneous ideal $I$ of the form $I(p_1)^{m_1}\cap\cdots\cap I(p_r)^{m_r}$,
where $p_1,\ldots,p_r\in\pr n$ are distinct points of \pr n and $m_1,\ldots,m_r$
are nonnegative integers (not all zero); we denote $Z$ by $m_1p_1+\cdots+m_rp_r$
and refer to $I(Z)=I$ as a {\it fat point\/} ideal, or as an {\it ideal of fat points}.
As in the case that $Z$ is smooth, the fat point ideal $I(Z)$
is defined by base point conditions;
i.e., the component $I(Z)_d$ of each degree $d$ of $I(Z)$
is the linear system of all forms of degree $d$ which vanish at each point $p_i$
with multiplicity at least $m_i$.

Easy examples show that a strict extension of the IGC does not hold for fat point
ideals; e.g., if $Z=2p_1+2p_2\subset \pr n$ with $n\ge 2$, then $I(Z)$ does
not have the maximal rank property, which raises the question
of how the number of generators in each degree can be determined for a fat point ideal. 
The ultimate goal of the line of research
which we initiate with this paper is to understand minimal homogeneous sets of generators
of ideals of the form $I(m_1p_1+\cdots+m_rp_r)$, where $m_1,\ldots,m_r$ are positive integers
and $p_1,\ldots,p_r$ are general points of \pr n. In the special case that $m_1=\cdots=m_r=m$,
we will refer to $Z=mp_1+\cdots+mp_r$ as a {\it uniform fat point subscheme}; the ideal 
$I(Z)$ is then the symbolic power $(I(p_1+\cdots+p_r))^{(m)}$.
In this paper we determine the number of generators of each degree in a minimal homogeneous 
set of generators for any symbolic power $(I(p_1+\cdots+p_r))^{(m)}$, where
$p_1,\ldots,p_r$ are $r\le9$ general points of \pr2, and we suggest a conjecture for $r>9$.

Since the Hilbert functions of fat point subschemes of \pr2 supported at 9 or fewer points
are known (see \ir{recall}, or, more generally, [\antican]), to determine numbers of
generators in each degree it is enough to determine the ranks of the maps $\mu_{d,1}$, so it is
this which will be our main concern. What we find (\ir{cokfact}) is that the cokernel of
$\mu_{d,1}$ has two parts. One part is related to fixed components, and its dimension can
be calculated by the same means as one computes the Hilbert functions of the fat point ideals
themselves, so it is the second part which is of most interest. 
For uniform fat point subschemes of \pr2 supported at $r\le 9$ general points,
we show that maximal rank holds for this second part except in three families of cases
(one family involving $r=7$ points, and two families involving $r=8$ points),
for which we explicitly determine the dimension of the second part (see \ir{subst}).

Taken together, our results allow one to recursively determine numbers of generators
in each degree in a minimal homogeneous set of generators of any symbolic power
$(I(p_1+\cdots+p_r))^{(m)}$ of an ideal $I(p_1+\cdots+p_r)$ corresponding to
$r\le9$ general points $p_1,\ldots,p_r$ of \pr2. Although results for $r>9$ are unknown,
we now discuss conjectures which would allow one to do the same in case $r>9$, and which
suggest that behavior for $r>9$ is even simpler.

In particular, let us say that the Generalized Ideal Generation Conjecture (GIGC) on \pr n holds
for $r$ if, for each $m>0$, the maximal rank property for $I(mp_1+\cdots+mp_r)$ holds
for general points $p_1,\ldots,p_r$ of \pr n. Then (as we will show) it turns out that:

\prclm{Theorem}{data}{For $r\le 9$, the GIGC on \pr 2
holds if and only if $r$ is 1, 4, or 9.}

The failure of the GIGC on \pr2 when $r$ is nonsquare less than 9 is by \ir{abfail}
guaranteed by the existence of
uniform abnormal curves for such $r$.
(Following Nagata [N1], a curve $C\subset\pr 2$ of degree $d$
whose multiplicity at each point $p_i$ is at least $m_i$ is said to be {\it abnormal\/} if 
$d\sqrt{r}<m_1+\cdots+m_r$, and {\it uniform\/}
if $m_1=\cdots=m_r$.) Nagata [N1] proves that
no abnormal curves occur for $r$ generic points when $r$ is a square, and he [N2] conjectures
that none occur for $r>9$. This prompts us, with some temerity perhaps, to 
propose the following conjecture:

\prclm{Conjecture}{VIGC}{The GIGC on \pr2 holds for all $r>9$.}

This also suggests the following question:

\prclm{Question}{ques}{Is there an $N$ depending on $n$, such that for each $r\ge N$ and each $m>0$,
$I(mp_1+\cdots+mp_r)$ has the maximal rank property for general points $p_1,\ldots,p_r$ of \pr n?}

In order to use \ir{VIGC} to actually determine numbers of generators for fat point subschemes
of \pr2, one needs to be able to determine their Hilbert functions.
Given $Z=mp_1+\cdots+mp_r\subset\pr2$ for $r>9$ sufficiently general points
$p_1,\ldots,p_r$, conjectures put forward in [\vanc] 
(equivalent to conjectures later put forward in [Hi], [Gi] and [\ravello]),
imply that $H_{I(Z)}(t)=\hbox{max}\{0,\binom{t+2}{2}-r\binom{m+1}{2}\}$
(i.e., each successive base condition $mp_i$ conjecturally imposes 
the expected number of additional independent conditions
on the remaining forms of degree $t$ as long as any forms remain). 
Via \ir{VIGC} above, one can now find $\nu_d$ for each $d$. 
Moreover, since we are working on \pr2, the numbers $\nu_d$ and the Hilbert function of $I(Z)$
determine the graded modules in a minimal free resolution
of $I(Z)$. For example, if $p_1,\ldots,p_{10}$ are general points of \pr2
and $Z=9(p_1+\cdots+p_{10})$, by the conjectures above there should be 15 generators
in degree 29, and 1 in degree 30. Denoting by $R$ the homogeneous coordinate ring
of \pr2, $I(Z)$ then would have the following minimal free resolution:

$$0\to R^{15}[-31]\to R^{1}[-30]\oplus R^{15}[-29]\to I(Z)\to 0.$$

To close this introduction, we discuss the significance of
$r=9$ as the boundary between what is understood and what is conjectural.
The approach taken in this paper is to work on the surface $X$ obtained by
blowing up the points $p_1,\ldots,p_r\in\pr2$, using cohomological methods 
applied to invertible sheaves. For any $r\le 9$ points the anticanonical class
$-K_X$ is the class of an effective divisor; in this situation the geometry
of divisors on $X$ is understood [\antican] and is used heavily in obtaining our 
results. For $r>9$ sufficiently general
points, $-K_X$ is not the class of an effective divisor and the geometry
of divisors on $X$ is not understood, raising a significant obstruction
to extending our results to $r>9$. In particular, whereas via [\antican] one can determine
the Hilbert function $H_{I(Z)}$ for any fat point subscheme $Z=m_1p_1+\cdots+m_rp_r\subset\pr2$
supported at any $r\le 9$ (even possibly infinitely near) points $p_1,\ldots,p_r$,
$H_{I(Z)}$ is unknown even for uniform fat point subschemes $Z$ supported at $r>9$ general points,
and determining numbers of generators is a more delicate question than that of determining $H_{I(Z)}$.

This paper is organized as follows. In \ir{surf} we accumulate
some notation and tools for working on surfaces. In \ir{main} we apply these tools
in our analysis of ideals of uniform fat point subschemes supported at
$r\le9$ general points of \pr2. Our analysis divides naturally into
three cases, $r\le 5$, $6\le r\le 8$, and $r=9$, which we treat separately. 

Hereafter, $R$ will denote the homogeneous
coordinate ring $k[\pr2]$ of \pr2.

\irrnSection{Background on Surfaces}{surf}
We will obtain our results on fat point ideals on \pr2 by working on
surfaces obtained by blowing up points of \pr2. We now establish the 
necessary connection.
Let $p_1,\ldots,p_r$ be distinct points of \pr2. Let $\pi:X\to \pr2$ 
be the morphism obtained by blowing up each point $p_i$. Let $E_i$ denote the exceptional
divisor of the blow up of $p_i$, and let $e_i$ denote its divisor class
(modulo linear equivalence). Let $e_0$ denote the pullback to $X$
of the class of a line in \pr2; the classes $e_0,\ldots,e_r$ comprise a ${\bf Z}$-basis
of $\hbox{Pic}(X)$ (where we identify divisor classes with their corresponding
invertible sheaves). Note that this basis is completely determined by $\pi$ and in turn
determines $\pi$. Also recall that 
$\hbox{Pic}(X)$ supports an intersection form, with respect to which
the basis $e_0,\ldots,e_r$ is orthogonal, satisfying
$-1=-e_0^2=e_1^2=\cdots=e_r^2$, and that the canonical class of $X$ is $K_X=-3e_0+e_1+\cdots+e_r$.

Consider now a fat point subscheme $Z=m_1p_1+\cdots+m_rp_r\subset\pr2$. 
Let $\C F_d$ denote the class
$de_0-m_1e_1-\cdots-m_re_r$. Since $e_0$ is the pullback
$\pi^*(\C O_{\pr2}(1))$ of the class of a line, 
we have for each $d$ and $i$ a natural isomorphism 
of $H^i(X,\C F_d)$ with 
$H^i(\pr2,\pi_*(-m_1e_1-\cdots-m_re_r)\otimes\C O_{\pr2}(d))=H^i(\pr2,\C I_Z(d))$.
In particular, the homogeneous coordinate ring 
$R=\bigoplus_{d\ge 0}H^0(\pr2,\C O_{\pr2}(d))$ can be identified
with $\bigoplus_{d\ge 0}H^0(X,de_0)$, and the
homogeneous ideal $I(Z)=\bigoplus_{d\ge 0}H^0(\pr2,\C I_Z(d))$ in $R$  can be identified with 
$\bigoplus_{d\ge 0}H^0(X,\C F_d)$. Moreover, under these identifications, the 
natural homomorphisms $H^0(X,\C F_d)\otimes H^0(X,e_0)\to H^0(X,\C F_{d+1})$
and $I(Z)_d\otimes R_1\to I(Z)_{d+1}$ correspond, so
the dimension $\nu_{d+1}$ of the cokernel of the latter is equal to
the dimension of the cokernel of the former. 

Following [Mu], we will
denote the kernel of $H^0(X,\C F_d)\otimes H^0(X,e_0)\to H^0(X,\C F_{d+1})$
by $\C R(\C F_d,e_0)$ and the cokernel by $\C S(\C F_d,e_0)$, and their dimensions
by $\r(\C F_d,e_0)$ and $\s(\C F_d,e_0)$.
Note that to say that $I(Z)_d\otimes R_1\to I(Z)_{d+1}$,
or equivalently that $H^0(X,\C F_d)\otimes H^0(X,e_0)\to H^0(X,\C F_{d+1})$,
has maximal rank is just to say that $[\r(\C F_d,e_0)][\s(\C F_d,e_0)]=0$.

Our main tool comes from [Mu]:

\prclm{Proposition}{Mumford}{Let $T$ be a closed subscheme of projective space,
let \C A and \C B be coherent sheaves on $T$ and let \C C be the class of an
effective divisor $C$ on $T$.
\item{(a)} If the restriction homomorphisms $H^0(T, \C A)\to H^0(C,\C A\otimes\C O_C)$
and $H^0(T, \C A\otimes\C B)\to H^0(C,\C A\otimes\C B\otimes\C O_C)$ are surjective
(for example, if $h^1(T,\C A\otimes\C C^{-1})=0=h^1(T,\C A\otimes\C C^{-1}\otimes\C B)$),
then we have an exact sequence
$$\eqalign{0\to & \C R(\C A\otimes\C C^{-1},\C B)\to\C R(\C A,\C B)\to\C R(\C A\otimes\C O_C,\C B)\to\cr
&\C S(\C A\otimes\C C^{-1},\C B)\to\C S(\C A,\C B)\to\C S(\C A\otimes\C O_C,\C B)\to0.\cr}$$
\item{(b)} If $H^0(T, \C B)\to H^0(C,\C B\otimes\C O_C)$ is surjective
(for example, if $h^1(T,\C B\otimes\C C^{-1})=0$), then 
$\C S(\C A\otimes\C O_C,\C B)=\C S(\C A\otimes\C O_C,\C B\otimes\C O_C)$.
\item{(c)} If $T$ is a smooth curve of genus $g$, and \C A and \C B are line 
bundles of degrees at least $2g+1$ and $2g$, respectively, then 
$\C S(\C A,\C B)=0$.}

\Prf See [Mu] for (a) and (c); we leave (b) as an easy exercise for the reader.\qed

Let \C F be the class of an effective divisor $F$ on a surface $X$. 
Let \C N denote the class of the fixed part of the linear system $|F|$;
then $\C H=\C F-\C N$ (called the {\it free part\/} of \C F) is fixed component free
and has $h^0(X, \C F)=h^0(X,\C H)$. 
The following lemma allows us to reduce a consideration of $\C S(\C F,e_0)$
to one of $\C S(\C H,e_0)$. 

\prclm{Lemma}{cokfact}{Let $e_0,\ldots,e_r$ be the divisor class group
basis corresponding to a blowing up $\pi:X\to\pr2$ at distinct points
$p_1,\ldots,p_r$, and let \C F be a divisor class on $X$. 
If \C F is not the class of an effective divisor, then 
$\s(\C F,e_0)=h^0(X,\C F+e_0)$.
If \C F is the class of
an effective divisor, let $\C H+\C N$ be its decomposition into free and fixed parts; 
then $\s(\C F,e_0)=[h^0(X,\C F+e_0)-h^0(X,\C H+e_0)]+\s(\C H,e_0)$.}

\Prf The case that \C F is not effective is clear, so assume that
\C F is the class of an effective divisor. 
Regarding $\C H$ and \C F as invertible sheaves, we have an inclusion
$\C H\to\C F$ which induces an isomorphism
on global sections. Thus we have a commutative diagram with exact columns:

$$\hbox{\hfil\vbox{\halign{\hfil $#$\hfil &\hfil $#$\hfil &\hfil $#$\hfil \cr
       0                          &     &           0            \cr
     \downarrow                   &     &      \downarrow        \cr
H^0(X,\C H)\otimes H^0(X,e_0) & \to & H^0(X,\C H+e_0)    \cr
     \downarrow                   &     &      \downarrow        \cr
H^0(X,\C F)\otimes H^0(X,e_0)     & \to & H^0(X,\C F+e_0)        \cr
     \downarrow                   &     &                        \cr
       0                          &     &                        \cr}}\hfil}$$

The image of $H^0(X,\C H)\otimes H^0(X,e_0) \to H^0(X,\C F+e_0)$ equals the image of 
$H^0(X,\C F)\otimes H^0(X,e_0)\to H^0(X,\C F+e_0)$ but factors through $H^0(X,\C H+e_0)$, 
which means that $\s(\C F,e_0)=
\s(\C H,e_0) + h^0(X,\C F+e_0)-h^0(X,\C H+e_0)$, giving the result. \qed

\rem{Remark}{recall} To determine generators for
the ideal $I(Z)$ of some uniform fat point subscheme $Z=m(p_1+\cdots+p_r)$ of \pr2,
it is enough by \ir{cokfact} on the blow up $X$ of \pr2 at $p_1,\ldots,p_r$ to
determine $h^0(X, de_0-m(e_1+\cdots+e_r))$ for every $d$, and, for each $d$ such that
$h^0(X, de_0-m(e_1+\cdots+e_r))>0$, to determine: the free part \C H of 
$de_0-m(e_1+\cdots+e_r)$; $\s(\C H,e_0)$; and $h^0(X, \C H+e_0)$.
(Since this also suffices to determine the Hilbert function of $I(Z)$,
it allows one in addition to write down a minimal free resolution of $I(Z)$,
as in the example near the end of \ir{intro} and in \ir{example}.)

In the case of any $r\le 9$ points, the results of [\antican] provide a solution
to determining $h^0(X,\C F)$ for any class \C F, and, when $h^0(X,\C F)>0$,
to finding the free part of \C F. For $r\le 9$ general points, these results 
are well known and can, in any case, be recovered from [\antican] or [\trans]; for the 
reader's convenience we recall the facts relevant to a uniform class
\C F in the case of $r$ general points of \pr2, first for $r\le 8$, and then for $r=9$.
(A class \C F on a blowing up $X$ of \pr2 at distinct points
$p_1,\ldots,p_r$ will be called a {\it uniform\/} class if
$\C F=de_0-m(e_1+\cdots+e_r)$ for some nonnegative integers $d$ and $m$.)

Let $X$ be the blowing up of $r\le 8$ general points of \pr2. If \C F is uniform and
if it is the class of an effective divisor, then the fixed part \C N is also uniform,
equal to $\Sigma t_{\C E}\C E$, where the sum is over all classes \C E of $(-1)$-curves and where
$t_{\C E}$ is the least nonnegative integer such that $(\C F-t_{\C E}\C E)\cdot \C E\ge 0$.
The classes of the $(-1)$-curves are known; up to permutation of the indices, 
they are (see Section 26 of [Ma]):
$e_1$, $e_0-e_1-e_2$, $2e_0-(e_1+\cdots+e_5)$, $3e_0-(2e_1+e_2+\cdots+e_7)$, 
$4e_0-(2e_1+2e_2+2e_3+e_4+\cdots+e_8)$, $5e_0-(2e_1+\cdots+2e_6+e_7+e_8)$, 
and $6e_0-(3e_1+2e_2+\cdots+2e_8)$. Now one can show that $\C N=0$ if $r=1$ or 4;
otherwise, \C N is
a nonnegative multiple of: $e_0-e_1-e_2$ if $r=2$; $3e_0-2e_1-2e_2-2e_3$
if $r=3$; $2e_0-(e_1+\cdots+e_5)$ for $r=5$; $12e_0-5(e_1+\cdots+e_6)$, $r=6$;
$21e_0-8(e_1+\cdots+e_7)$, $r=7$; or $48e_0-17(e_1+\cdots+e_8)$, $r=8$.
It also follows that a uniform class $de_0-m(e_1+\cdots+e_r)$ is the class of an effective
divisor if and only if $d\ge \epsilon_r m$, where $\epsilon_1=\epsilon_2=1$,
$\epsilon_3=3/2$, $\epsilon_4=\epsilon_5=2$, $\epsilon_6=12/5$, $\epsilon_7=21/8$,
and $\epsilon_8=48/17$.

Recall that a class being {\it numerically effective\/} means that in the intersection form 
it meets every effective divisor nonnegatively.
In particular, the free part of the class of an effective divisor 
is numerically effective; conversely, if $X$ is any
blowing up of \pr2 at 8 or fewer points, general or not [\antican], then
any numerically effective class
\C F on $X$ is the class of an effective divisor with no fixed components 
and has $h^1(X,\C F)=h^2(X,\C F)=0$, hence $h^0(X,\C F)=(\C F^2-K_X\cdot \C F)/2+1$
by Riemann-Roch for surfaces.

Finally, let $X$ be the blowing up of $r=9$ general points of \pr2. 
Then there is a unique smooth cubic through the 9 points, so
$-K_X=3e_0-e_1-\cdots-e_9$ is numerically effective, the class of a smooth elliptic curve.
If \C F is uniform, 
we can write $\C F=te_0-sK_X$ for uniquely determined
integers $t$ and $s\ge 0$, and \C F is the class of an effective divisor
if and only if $t$ is also nonnegative, in which case $h^1(X,\C F)=h^2(X,\C F)=0$,
hence $h^0(X,\C F)=(\C F^2-K_X\cdot \C F)/2+1$. Moreover, given $s\ge 0$, if $t>0$, then
\C F is fixed part free, while if $t=0$, then $h^0(X,\C F)=1$.
\qed

The next result will be helpful in verifying failure of the GIGC.
Call a uniform class $\C E=de_0-m(e_1+\cdots+e_r)$ 
on a blowing up $X$ of \pr2 at distinct points
$p_1,\ldots,p_r$ {\it abnormal\/} if \C E is the class of 
an effective divisor with $d<\sqrt{r}m$
(note that this is equivalent to $\C E^2<0$).

\prclm{Proposition}{abfail}{Let $X$ be a blowing up
of $r$ distinct points $p_1,\ldots,p_r$ of \pr2. If $X$ has
a uniform abnormal class \C E, then for some positive integers
$\alpha$ and $\beta$, $I(\beta(p_1+\cdots+p_r))_\alpha
\otimes R_1\to I(\beta(p_1+\cdots+p_r))_{\alpha+1}$
does not have maximal rank.}

\Prf Since \C E is the class of an effective divisor of negative self-intersection,
we can find positive integers $a$ and $b$ such that $ae_0+b\C E$ has nontrivial fixed part
but such that $(a+1)e_0+b\C E$ has trivial fixed part. Now,
$ae_0+b\C E=\alpha e_0-\beta(e_1+\cdots+e_r)$ for some positive $\alpha$ and $\beta$.
Since $a>0$, $H^0(X,ae_0)\otimes H^0(X,e_0)\to H^0(X,(a+1)e_0)$ is not injective,
hence neither is $H^0(X,ae_0+b\C E)\otimes H^0(X,e_0)\to H^0(X,(a+1)e_0+b\C E)$.
Since $(a+1)e_0+b\C E$ is fixed component free but $ae_0+b\C E$ is not,
we see $H^0(X,ae_0+b\C E)\otimes H^0(X,e_0)\to H^0(X,(a+1)e_0+b\C E)$
is also not surjective. Thus $H^0(X,ae_0+b\C E)\otimes H^0(X,e_0)\to H^0(X,(a+1)e_0+b\C E)$,
and hence $I(\beta(p_1+\cdots+p_r))_\alpha
\otimes R_1\to I(\beta(p_1+\cdots+p_r))_{\alpha+1}$, do not have maximal rank.\qed

The following result is well known (see Proposition 3.7 of [DGM]) and follows easily
by appropriately applying \ir{Mumford}.

\prclm{Lemma}{myomega}{Let $e_0,\ldots,e_r$ be the classes
corresponding to a blowing up $X\to\pr2$ at
distinct points $p_1,\ldots, p_r$. Let $Z=m_1p_1+\cdots+m_rp_r$,
and let $\C F_d$ denote $de_0-m_1e_1-\cdots-m_re_r$.
If $\omega_Z$ is the degree of a generator of greatest degree in a minimal
homogeneous set of generators of $I(Z)$ (equivalently, $\mu_d$ is surjective
for $d\ge \omega_Z$ but not for $d=\omega_Z-1$) and if $\tau_Z$ is the least integer
such that $h^1(X,\C F_t)=0$ for $t\ge\tau_Z$, then  $\omega_Z\le \tau_Z+1$.
In particular, $\C S(\C F_t,e_0)=0$ for $t>\tau_Z$.}

Indeed, $\tau_Z+1$ is just the regularity of $I(Z)$.

\irrnSection{Main Results}{main}
We now work out our results for ideals of uniform fat point subschemes supported at $r\le 9$
general points of \pr2. We divide our analysis into three cases, $r\le 5$, $6\le r\le 8$,
and $r=9$, with the second case requiring most of the effort but also being the most interesting.

\irSubsection{Five or Fewer General Points}{conic}
By \ir{cokfact} and \ir{recall}, we are reduced to determining
$\C S(\C F,e_0)$ for numerically effective classes \C F on
the blow up $X$ of \pr2 at $r\le 5$ general points.
But any five or fewer general points in the plane lie on a smooth conic,
so the results of [Cat] apply. Translating the results of [Cat] to the language used here and
examining what [Cat] proves, we find that 
$\C S(\C F,e_0)=0$ for any numerically effective class \C F. 
(In fact, [Cat] iteratively finds generators for and a resolution of
$I(Z)$ for any fat point subscheme $Z=m_1p_1+\cdots+m_tp_t$, where
$p_1,\ldots,p_t$ are distinct points of a smooth plane conic,
which includes the case of a uniform $Z$ supported at 5 or fewer general points
of \pr2. From our perspective, the key fact in [Cat], not explicitly stated there, is that
$\C S(\C F,e_0)=0$ for any numerically effective class \C F on the blow up 
$X$ of points on a smooth conic, which reduces one, by \ir{cokfact}, to cohomology
calculations already carried out in [\vanc]. Moreover, the key fact generalizes;
see \ir{CatGen} below.)

Applying the foregoing to $Z=m(p_1+\cdots+p_r)$ for
$r\le 5$ general points $p_1,\ldots,p_r\in\pr2$ and $m>0$, we have the following.
For $r=1$, it is easy to see that $I(Z)_t=0$ for $t<m$ and that $te_0-me_1$ is
numerically effective for $t\ge m$. 
The former means that $\mu_t$ is injective for $t<m$, and by the preceding paragraph
and numerical effectivity of $te_0-me_1$ for $t\ge m$, we have 
$\C S(te_0-me_1,e_0)=0$ for $t\ge m$, and hence $\mu_t$ is surjective for $t\ge m$. 
Thus the GIGC holds on \pr2 for $r=1$.
For $r=4$, $I(Z)_t=0$ for $t<2m$, since $2e_0-(e_1+\cdots+e_4)$ is 
numerically effective but $[2e_0-(e_1+\cdots+e_4)]\cdot[te_0-m(e_1+\cdots+e_4)]<0$.
Also, $\C S(\C F_t,e_0)=0$ for $\C F_t=te_0-m(e_1+\cdots+e_4)$ and $t\ge 2m$,
since $\C F_t=m(2e_0-(e_1+\cdots+e_4))+(t-2m)e_0$ is numerically effective. 
Thus the GIGC holds on \pr2 also for $r=4$.  

To see that the GIGC on \pr2 fails for $r=2,3,5$, it is enough by \ir{abfail}
to find in each case a uniform abnormal class. But these have already been exhibited in
\ir{recall}: for $r=2$, we have
$e_0-(e_1+e_2)$; for $r=3$, there is $3e_0-2(e_1+e_2+e_3)$;
and for $r=5$, $2e_0-(e_1+\cdots+e_5)$. One can check, in fact, that for $r=2,3,5$,
$I(m(p_1+\cdots+p_r))$ fails to have the maximal rank property if and only if:
$r=2$ and $m\ge 2$; or $r=3$ or $r=5$ and $m\ge 3$.

\rem{Remark}{CatGen} It turns out that the result that $\C S(\C F,e_0)=0$
for any numerically effective class on a blow up $X$ of points on a smooth
conic is true more generally. In fact, let $X\to\pr2$ be any projective birational
morphism where $X$ is smooth and projective (i.e., $X$ is obtained by blowing up points,
possibly infinitely near): If $-K_X-e_0$ is the class of an effective divisor
(this is always true for a blowing up of points on a smooth conic),
then $\C S(\C F,e_0)=0$ for any numerically effective class \C F on $X$ [\fatpts].
So, by \ir{cokfact}, for any divisor class \C G on such an $X$ we have 
$\s(\C G,e_0)=h^0(X,\C G+e_0)$ if \C G is not the class of an effective divisor,
or, if it is, $\s(\C G,e_0)=h^0(X,\C G+e_0)-h^0(X,\C G'+e_0)$,
where $\C G'$ is the free part of \C G; and using
the results of [\antican], one can compute $h^0$ and the free part
for any divisor class on $X$. [The condition that 
$-K_X-e_0$ be the class of an effective divisor just means that, in an appropriate
sense, $X$ is obtained by blowing up points on a conic (but the conic need not be
irreducible and the points can be infinitely near). Thus for points $p_1,\ldots,p_r$, possibly
infinitely near, of a plane conic, possibly reducible or nonreduced, one can explicitly find 
a minimal set of generators for $I(Z)$, the Hilbert function of $I(Z)$, and
a minimal free resolution for $I(Z)$, where $Z$ is 
any fat point subscheme $Z=m_1p_1+\cdots+m_rp_r$ (see [\fatpts]).]

\irSubsection{Six to Eight General Points}{sixtoeight}
We begin by pointing out that the GIGC fails by \ir{abfail} for $6\le r\le 8$, since,
as shown in \ir{recall}, in each case
the blow up of $r$ general points supports a uniform abnormal class:
for $r=6$, $\C E=12e_0-5(e_1+\cdots+e_6)$ is such;
for $r=7$, $\C E=21e_0-8(e_1+\cdots+e_7)$ is such; and 
for $r=8$, $\C E=48e_0-17(e_1+\cdots+e_8)$ is such. 

Thus the maximal rank property need not hold for ideals of the form
$I(Z)$ where $Z=m(p_1+\cdots+p_r)$ and $p_1,\ldots,p_r$ are $6\le r\le 8$ general points
of \pr2, so, as discussed above, in order to find the number of generators in each degree
in a minimal homogeneous set of generators, we need only find $\s(\C F, e_0)$ for numerically
effective uniform classes \C F on the blow up $X$ of the points. Unlike what happens for
$r\le 5$, however, $\s(\C F, e_0)$ need not vanish; we have instead:

\prclm{Theorem}{subst}{Let $\C F=\C F(d,m,r)$ be a uniform numerically effective class
on the blowing up of $6\le r\le 8$ general points of \pr2 (where $\C F(d,m,r)$
denotes $de_0-m(e_1+\cdots+e_r)$). 
\itemitem{(a)} If $r=6$, then $\r(e_0,\C F)\s(e_0,\C F)=0$.
\itemitem{(b)} If $r=7$, then $\r(e_0,\C F)\s(e_0,\C F)=0$ unless 
$\C F=l\C F(8,3,7)$ for $l\ge 3$, in which case $\s(e_0,\C F)=7$.
\itemitem{(c)} If $r=8$, then $\r(e_0,\C F)\s(e_0,\C F)=0$, unless
$\C F=l\C F(17,6,8)$ for $l\ge 9$, in which case $\s(e_0,\C F)=48$, or
unless $\C F=l\C F(17,6,8)+\C F(3,1,8)$ for 
$l\ge 6$, in which case $\s(e_0,\C F)=16$.}

\Prf Let $6\le r\le 8$ and let $\C F=de_0-m(e_1+\cdots+e_r)$ be a uniform
class. If \C F is numerically effective, then $h^1(X,\C F+te_0)=0$ 
for all $t\ge 0$ (by \ir{recall}), so $\C S(\C F+te_0,e_0)=0$
for all $t>0$ by \ir{myomega}. Thus we only need to consider $\C F=\delta e_0-m(e_1+\cdots+e_r)$,
where $\delta$ is the least $d$ such that $de_0-m(e_1+\cdots+e_r)$ is numerically effective.
Using \ir{recall} it follows that $\delta$
is the least positive integer $d$ such that: $d\ge 12m/5$ if $r=6$;
$d\ge 21m/8$ if $r=7$; or $d\ge 48m/17$ if $r=8$. 

First say $r=6$ and let $\C F=\delta e_0-m(e_1+\cdots+e_r)$.
If $m$ is odd, then $\C F=-K_X+(m-1)(5e_0-2(e_1+\cdots+e_6)/2$,
while $\C F=m(5e_0-2(e_1+\cdots+e_6))/2$ if $m$ is even. 
In any case $h^2(X,\C F-e_0)=0$ by duality.

If $m$ is odd, one checks that $h^1(X,\C F-e_0)=0$,
and hence (by \ir{myomega}) that $\C S(\C F,e_0)=0$.
(Explicitly, if $m$ is odd, then each of the six
$(-1)$-curves $2e_0-(e_1+\cdots+e_6)+e_i$ meets $\C F-e_0$
negatively; thus either we have $\C E=12e_0-5(e_1+\cdots+e_6)$
in the fixed part of $\C F-e_0$ or we have $h^0(X,\C F-e_0)=0$.
The latter happens for $m=1,3$, and the former if $m\ge 5$,
in which case it is easy to check that $\C F-e_0-\C E$ is numerically effective
and hence that $h^1(X, \C F-e_0-\C E)=0$. For $m=1,3$,
$h^1(X,\C F-e_0)=0$ now follows since $h^0-h^1+h^2$ applied to
$\C F-e_0$ vanishes by Riemann-Roch; $h^1(X,\C F-e_0)=0$ also follows
for $m\ge 5$, since $h^0(X,\C F-e_0-\C E)=h^0(X,\C F-e_0)$
and by Riemann-Roch and arithmetic $h^0-h^1+h^2$ gives the same result
applied either to $\C F-e_0-\C E$ or to $\C F-e_0$.)

Now suppose $m=2s$, with $s\ge 2$. For $s=2$,
$e_0\cdot(\C F-\C E+e_0)=-1$, so $h^0(X,\C F-\C E+e_0)=0$
so $\C S(\C F-\C E,e_0)=0$. For $s>2$, $\C F-\C E$
is numerically effective with odd uniform multiplicity,
so by the preceding paragraph $\C S(\C F-\C E,e_0)=0$.
Let $E$ be the effective divisor whose class
is \C E. Then $\C F\otimes\C O_E=\C O_E$,
and using \ir{Mumford}(b) it is easy to check that 
$\C S(\C O_E,e_0)=0$. If we check that $h^1(X,\C F-\C E+e_0)=0$ and
$h^1(X,\C F-\C E)=0$, then
we can apply \ir{Mumford}(a)
to $(0\to\Gamma(\C F-\C E)\to\Gamma(\C F)\to\Gamma(\C F\otimes\C O_E)\to0)\otimes\Gamma(e_0)$
to obtain $\C S(\C F,e_0)=0$.
But for $s>2$, we have $h^1(X,\C F-\C E+e_0)=0$ and
$h^1(X,\C F-\C E)=0$ by \ir{recall}. For $s=2$, we have $\C F-\C E=K_X+e_0$
and $\C F-\C E+e_0=K_X+2e_0$, so using duality and descending to \pr2 we see
$h^1(X,\C F-\C E+ae_0)=h^1(\pr2,\C O_{\pr2}(1-a))=0$ for any $a$.

We are left with the case $m=2$; apply \ir{Mumford} 
to $(0\to\Gamma(\C F-\C C)\to\Gamma(\C F)\to\Gamma(\C F\otimes\C O_C)\to0)\otimes\Gamma(e_0)$,
where $C$ is the $(-1)$-curve whose class is $\C C=2e_0-e_1-\cdots-e_5$.
Since $h^0(X, \C F-\C C)=2$ and $e_0\cdot(\C F-\C C)=3$, the sections of $\C F-\C C$
correspond to a pencil of cubic plane curves, so
we see $\r(\C F-\C C,e_0)=0$ (i.e., any nontrivial element of the kernel of 
$H^0(X, \C F-\C C)\otimes H^0(X,e_0)\to H^0(X, \C F-\C C+e_0)$
must be of the form $f\otimes l_1-g\otimes l_2$, where $f$ and $g$ define distinct
elements of the cubic pencil and $l_1$ and $l_2$ define distinct lines,
which means that $f$ and $g$ have a factor in common corresponding to
a plane conic curve, but $\C F-\C C$ is numerically effective so
fixed part free by \ir{recall}, so no such common factor is possible), so 
$\s(\C F-\C C,e_0)=h^0(X, \C F-\C C+e_0)-h^0(X, \C F-\C C)h^0(X,e_0)$
which we compute to be 1. By a similar dimension count,
we see $\s(\C F,e_0)\ge 1$, but by \ir{Mumford}(b),(c) we
have $\s(\C F\otimes\C O_C,e_0)=0$. Thus  $\C S(\C F-\C C,e_0)$ surjects
onto $\C S(\C F,e_0)$ by \ir{Mumford}(a) so in fact $\s(\C F,e_0)=1$,
hence $\r(\C F,e_0)=0$, whence $\r(\C F,e_0)\s(\C F,e_0)=0$ as claimed.

The cases of $r=7$ and 8 general points of \pr 2
are handled similarly so we indicate these in less detail.
So now say $r=7$; given an integer $m>0$,
let $\C F_m$ denote $\delta e_0-m(e_1+\cdots+e_r)$,
where $\delta$ is the least positive integer $d$ such that
$de_0-m(e_1+\cdots+e_r)$ is numerically effective but $(d-1)e_0-m(e_1+\cdots+e_r)$ is not.
Write $m=3s+t$, where $0\le t<3$ is the remainder when $m$ is divided by 3.
Then $\C F_m=s\C F_3-tK_X$, where $\C F_3=8e_0-3(e_1+\cdots+e_7)$.
Also, let $E$ be the effective divisor whose class is 
$\C E=21e_0-8(e_1+\cdots+e_7)$; $E$ is a union of seven disjoint $(-1)$-curves. 

For $t=2$, we have $h^1(X, \C F_m-e_0)=0$ which gives $\C S(\C F_m,e_0)=0$, 
so now consider $t=1$. If $s\ge 2$, then $h^1(X,\C F_m-\C E+e_0)=h^1(X,\C F_m-\C E)=0$,
so, if we show that $\C S(\C F_m-\C E,e_0)$ and $\C S(\C F_m\otimes\C E,e_0)$ 
vanish, we conclude $\C S(\C F_m,e_0)=0$ by applying \ir{Mumford}(a) to
$(0\to\Gamma(\C F_m-\C E)\to\Gamma(\C F_m)\to\Gamma(\C F_m\otimes\C O_E)\to0)\otimes\Gamma(e_0)$.
For $s\ge 3$ we have $\C F_m-\C E=\C F_{3(s-3)+2}$, for which we have already established
$\C S(\C F_m-\C E,e_0)=0$, while for $s=2$ we have $h^0(X,\C F_m-\C E+e_0)=0$,
so again $\C S(\C F_m-\C E,e_0)=0$. Also, $\C S(\C F_m\otimes\C O_E,e_0)=0$:
$E$ is a disjoint union of seven $(-1)$-curves, and for each of these curves $C$
we have $\C F_m\otimes\C O_C=O_C(1)$, so it is enough to show
that $\Gamma(O_C(1))\otimes\Gamma(e_0)\to\Gamma(O_C(1)\otimes e_0)=\Gamma(O_C(4))$
is surjective. But the restriction $V\subset \Gamma(O_C(3))$ 
of $\Gamma(e_0)$ to $C$ is a base point free cubic web,
and it is not hard to check for such a $V$ that $\Gamma(O_C(1))\otimes V$ surjects
to $\Gamma(O_C(4))$, as required.

In the $t=1$ case we are left with treating $s=1$ and $s=0$. 
For $s=1$, we have $m=4$,
which we handle by applying \ir{Mumford}(a) to 
$(0\to\Gamma(2e_0-e_4-\cdots-e_7)\to\Gamma(\C F_4)\to
\Gamma(\C F_4\otimes\C O_D)\to0)\otimes\Gamma(e_0)$,
where $D$ is a section of $9e_0-4e_1-4e_2-4e_3-3e_4-\cdots-3e_7$.
The conclusion now follows since $\C S(2e_0-e_4-\cdots-e_7,e_0)=0$
by \ir{conic}, and $\C S(\C F_4\otimes\C O_D,e_0)=0$
by an argument similar to that used above for $s\ge 2$, keeping in mind
that $D$ is a union of three of the seven $(-1)$-curves comprising $E$.
For $s=0$, we have $m=1$,
so $\C F_m=-K_X$, which we handle by applying \ir{Mumford}(a),(c) to 
$(0\to\Gamma(\C O_X)\to\Gamma(\C F_m)\to\Gamma(\C F_m\otimes\C O_D)\to0)\otimes\Gamma(e_0)$,
where now $D$ is a smooth section of $-K_X$ (and hence a smooth elliptic curve).

There remains only the $t=0$ case. For $s\ge 3$, we have
$\C F_{3s}-\C E=\C F_{3(s-3)+1}$, so 
$\C S(\C F_{3s}-\C E,e_0)=0$.
Also, $\s(\C F_{3s}\otimes\C O_E,e_0)=7$, since $E$ has seven components
and on each component $C$ we have $\C F_{3s}\otimes\C O_C=\C O_C$,
so the images of $H^0(C,\C O_C)\otimes H^0(X,e_0)\to H^0(C,\C O_C(3))$
and $H^0(X,e_0)\to H^0(C,\C O_C(3))$ are equal, but the latter
has 1 dimensional cokernel. Now we apply \ir{Mumford} to
$(0\to\Gamma(\C F_m-\C E)\to\Gamma(\C F_m)\to\Gamma(\C F_m\otimes\C O_E)\to0)\otimes\Gamma(e_0)$
to conclude $\s(\C F_{3s},e_0)=7$. 

We are left with $s=1$ and 2, which we handle by applying \ir{Mumford} to
$(0\to\Gamma(\C F_{3s}-\C D)\to\Gamma(\C F_{3s})\to
\Gamma(\C F_{3s}\otimes\C O_D)\to0)\otimes\Gamma(e_0)$,
using $\C D=9e_0-4e_1-4e_2-4e_3-3e_4-\cdots-3e_7$ (as above)
for $s=1$ (obtaining $\r(\C F_{3}-\C D,e_0)=0=\r(\C F_{3}\otimes\C O_D,e_0)$
and hence $\r(\C F_{3},e_0)=0$, from which we compute $\s(\C F_{3},e_0)=4$),
and using $\C D=18e_0-7e_1-\cdots-7e_6-6e_7$ for $s=2$
(obtaining $\r(\C F_{6}-\C D,e_0)=0=\r(\C F_{6}\otimes\C O_D,e_0)$
and hence $\r(\C F_{6},e_0)=0$ and so $\s(\C F_{6},e_0)=6$).

In conclusion, for $r=7$, $\r(\C F,e_0)\s(\C F,e_0)=0$ for all
numerically effective uniform classes \C F except $\C F=l\C F_3$
for $l\ge3$, in which case $\s(\C F,e_0)=7$.

We now proceed to the last case, for which $X$ is a blow up of
$r=8$ general points of \pr2. Here we let $E$ be the effective divisor whose class is 
$\C E=48e_0-17(e_1+\cdots+e_8)$; $E$ is a union of eight disjoint $(-1)$-curves, each of
which comes from a plane sextic with seven double points and a triple point.
Using notation analogous to that above, $\C F_m$ denotes the class
$\delta e_0-m(e_1+\cdots+e_8)$, where $\C F_m$ is numerically effective but
$\C F_m-e_0$ is not. Similarly to what was done above, we can write 
$\C F_m=s\C F_6-tK_8$, where $\C F_6=17e_0-6(e_1+\cdots+e_8)$ and 
$0\le t<6$ is the remainder when $m$ is divided by $s$.

Following the pattern for 7 points, we find $h^1(X,\C F_{6s+5}-e_0)=0$ for $s\ge 0$, hence
$\C S(\C F_{6s+5},e_0)=0$ for all $s\ge 0$, and we find $\C F_{6s+t}-\C E=\C F_{6(s-3)+t+1}$
and hence $h^1(X, \C F_{6s+t}-\C E+e_0)=h^1(X, \C F_{6s+t}-\C E)=0$ if $s\ge3$ and $0\le t<5$.
Thus for $m\ge18$ we can apply \ir{Mumford} to 
$(0\to\Gamma(\C F_m-\C E)\to\Gamma(\C F_m)\to\Gamma(\C F_m\otimes\C O_E)\to0)\otimes\Gamma(e_0)$.
To do so we will need to determine $\s(\C F_m\otimes\C O_E,e_0)$
(or, equivalently, $\r(\C F_m\otimes\C O_E,e_0)$).
As we might expect from the case of seven points,
if \C C is the class of any component $C$ among the eight components of $E$, then
$\s(\C F_m\otimes\C O_E,e_0)=8\s(\C F_m\otimes\C O_C,e_0)$, so we may restrict our attention
to $C$. Note that $\C F_m\otimes\C O_C=-tK_8\otimes\C O_C=\C O_C(t)$.
For $t=0$, clearly $\r(\C O_C,e_0)=0$ (a linear form times a nonzero
constant cannot vanish on a sextic) whence $\s(\C O_C,e_0)=4$. 

For $t=1$, again $\r(\C O_C(1),e_0)=0$, whence $\s(\C O_C(1),e_0)=2$:
letting $x$ and $y$ be a basis for $\Gamma(\C O_C(1))$, 
a nontrivial element of $\Gamma(\C O_C(1))\otimes\Gamma(e_0)$ which maps to 0 in 
$\Gamma(\C O_C(1)\otimes e_0)=\Gamma(\C O_C(7))$ gives an equation
$xf=yg$, where $f$ and $g$ are restrictions to $C$ of distinct lines in \pr2.
But $f$ and $g$ have degree 6, so $xf=yg$ implies $f$ and $g$ have 5 zeros on $C$ in common.
Since the image of $C$ in \pr2 has at most a triple point, two distinct
lines can have at most 3 points of $C$ in common, contradicting there being a nontrivial
element of the kernel.

For $t=2$, both $\C R(\C O_C(2),e_0)$ and $\C S(\C O_C(2),e_0)$ vanish:
let $x$ and $y$ be as before and let $f,g,h$ be a basis for the restriction
of $\Gamma(e_0)$ to $C$ such that $f$ and $g$ correspond to lines in \pr2
which meet at the triple point of the image of $C$ in \pr2.
If $\r(\C O_C(2),e_0)\ne0$, then we have an equation
$q_1f+q_2g+q_3h=0$, where $q_1,q_2,q_3$ (not all 0) lie in the span of
$\{x^2,xy,y^2\}$. Since $f$ and $g$ have exactly 3 zeros in common,
we cannot have $q_3=0$, and so $h$ also has a zero in common with $f$ and $g$,
which gives the contradiction that the restriction of $\Gamma(e_0)$ to $C$
has a base point. Thus $\r(\C O_C(2),e_0)=0$ from which we easily compute
$\s(\C O_C(2),e_0)=0$.

For $t=3,4$ or 5, we have $\s(\C O_C(t),e_0)=0$: say $t=3$ ($t=4$ or 5 are similar).
Let $x$ and $y$ be as above; thus cubics in $x$ and $y$ span $\Gamma(\C O_C(3))$.
But $\Gamma(\C O_C(1))\otimes\Gamma(\C O_C(2))$ surjects onto
$\Gamma(\C O_C(3))$, and, by the previous case,
$\Gamma(\C O_C(2))\otimes\Gamma(e_0)$ surjects onto
$\Gamma(\C O_C(8))$, so $\Gamma(\C O_C(3))\otimes\Gamma(e_0)$
and $(\Gamma(\C O_C(1))\otimes\Gamma(\C O_C(2)))\otimes\Gamma(e_0)$
and $\Gamma(\C O_C(1))\otimes\Gamma(\C O_C(8))$
all have the same image in $\Gamma(\C O_C(9))$. Since $C$ is rational,
we know $\Gamma(\C O_C(1))\otimes\Gamma(\C O_C(8))$
surjects onto $\Gamma(\C O_C(9))$, whence $\s(\C O_C(t),e_0)=0$.

Now apply \ir{Mumford} to 
$(0\to\Gamma(\C F_m-\C E)\to\Gamma(\C F_m)\to\Gamma(\C F_m\otimes\C O_E)\to0)\otimes\Gamma(e_0)$.
For $m=6s+4$, $\s(\C F_m\otimes\C O_E,e_0)=0$ and for $s\ge 3$ 
we have $\C F_m-\C E=\C F_{6(s-3)+5}$, so $\s(\C F_m-\C E,e_0)=0$; thus 
for $t=4$ and $s\ge 3$ we have $\s(\C F_m,e_0)=0$. To handle $s=0$,
note that $\C F_4=-4K_X$, and indeed, for $m<6$, $\C F_m=-mK_X$.
But for $-mK_X$ with $m<5$ we apply \ir{Mumford}(a,b,c) to
$(0\to\Gamma(-(m-1)K_X)\to\Gamma(-mK_X)\to\Gamma(-mK_X\otimes\C O_C)\to0)\otimes\Gamma(e_0)$,
where $C$ is a smooth section of $-K_X$, hence an elliptic curve.
It follows that $\s(-K_X,e_0)=1$ and $\r(-K_X,e_0)=0$, and, if $\s(-2K_X,e_0)=0$, that $\s(-mK_X,e_0)=0$
for $m=3$ and $4$ (we already know $\s(-mK_X,e_0)=0$ for $m=5$). To check
$\s(-2K_X,e_0)=0$, apply \ir{Mumford} to
$(0\to\Gamma(-2K_X-\C C)\to\Gamma(-2K_X)\to\Gamma(-2K_X\otimes\C O_C)\to0)\otimes\Gamma(e_0)$,
where $C$ is the $(-1)$-curve whose class is
$6e_0-3e_1-2e_2-\cdots-2e_8$.

To handle $t=4$ and $s=1$, apply \ir{Mumford} to
$(0\to\Gamma(\C C+K_X)\to\Gamma(\C C)\to\Gamma(\C C\otimes\C O_D)\to0)\otimes\Gamma(e_0)$,
where $D$ is a smooth section of $-K_X$, $\C C=\C F_{10}-\C G$ and $\C G=
24e_0-9(e_1+\cdots+e_4)-8(e_5+\cdots+e_8)$ is the class of the union of four
disjoint sextic $(-1)$-curves. This shows $\s(\C C,e_0)=0$, which by considering
$(0\to\Gamma(\C F_{10}-\C G)\to\Gamma(\C F_{10})\to\Gamma(\C F_{10}\otimes\C O_G)\to0)\otimes\Gamma(e_0)$,
where $G$ is the effective divisor whose class is \C G, shows $\s(\C F_{10},e_0)=0$.
The case of $t=4$ and $s=2$ is similar, but with $\C G=42e_0-15(e_1+\cdots+e_7)-14e_8$
the union of seven $(-1)$-curves. 

In a similar way, since now $\s(\C F_{6s+4},e_0)=0$ for all $s\ge 0$,
we find $\s(\C F_{6s+3},e_0)=0$ for all $s\ge 3$ and we are left with the cases
$s=1$ and 2 ($s=0$ was done above). For $s=1$, consider
$(0\to\Gamma(\C F_m-\C C)\to\Gamma(\C F_m)\to\Gamma(\C F_m\otimes\C O_C)\to0)\otimes\Gamma(e_0)$
with $\C C=24e_0-9(e_1+\cdots+e_4)-8(e_5+\cdots+e_8)$; for $s=2$, replace
$24e_0-9(e_1+\cdots+e_4)-8(e_5+\cdots+e_8)$ by $42e_0-15(e_1+\cdots+e_7)-14e_8$.
We find for all $s\ge 0$ that $\s(\C F_{6s+3},e_0)=0$.

With this in hand, following the same pattern for $t=2$ we find $\s(\C F_{6s+2},e_0)=0$ for all $s\ge 3$.
We have done the $s=0$ case above; for $s=1$ consider first 
$(0\to\Gamma(\C F_8-\C C-\C D)\to\Gamma(\C F_8-\C C)\to\Gamma((\C F_8-\C C)\otimes\C O_D)\to0)\otimes\Gamma(e_0)$
where $\C D=2e_0-e_4-\cdots-e_8$ is a $(-1)$-curve and
$\C C=18e_0-7(e_1+e_2+e_3)-6(e_4+\cdots+e_8)$ is the union of three $(-1)$-curves, then consider
$(0\to\Gamma(\C F_8-\C C)\to\Gamma(\C F_8)\to\Gamma(\C F_8\otimes\C O_C)\to0)\otimes\Gamma(e_0)$, 
to show $\r(\C F_8,e_0)=0$ and $\s(\C F_8,e_0)=1$. 

Now consider $s=2$; i.e., $\C F_{14}$.
Using $(0\to\Gamma(\C F_{14}-\C C)\to\Gamma(\C F_{14})\to\Gamma(\C F_{14}\otimes\C O_C)\to0)\otimes\Gamma(e_0)$
with $\C C=36e_0-13(e_1+\cdots+e_6)-12(e_7+e_8)$, we find by \ir{Mumford} that
$\s(\C G,e_0)=\s(\C F_{14},e_0)$, where $\C G=\C F_{14}-\C C=4e_0-(e_1+\cdots+e_6)-2(e_7+e_8)$.
To see that $\s(\C G,e_0)=0$, note that $H^0(X,\C G)\otimes H^0(X,e_0)$ and
$H^0(X, \C G+e_0)$ both have dimension 9. Let $V\subset H^0(X, \C G)$ be the linear
subsystem with a fixed component corresponding to a line through $p_7$ and $p_8$;
i.e., $V$ is the cubic pencil of sections of $-K_X$ plus a fixed component whose class is $e_0-e_7-e_8$.
Note that the image $V'$ of $V\otimes H^0(X, e_0)$ in $H^0(X, \C G+e_0)$ has dimension 6 and a base point
off the fixed component (the base point being that of the cubic pencil). Let $W\in H^0(X, \C G)$ be the element
associated to the divisor whose two components are the $(-1)$-curves whose classes
are $2e_0-e_1-e_2-e_3-e_7-e_8$ and $2e_0-e_4-e_5-e_6-e_7-e_8$.
Let $W'$ be the image of $W\otimes H^0(X, e_0)$ in $H^0(X,\C G+e_0)$; by dimension count, if 
$V'$ and $W'$ have only 0 in common, then $\s(\C G,e_0)=0$. But any element of $W'$
corresponds to a conic $Q_1$ through $p_1,p_2,p_3,p_7,p_8$, another conic $Q_2$ through $p_4,p_5,p_6,p_7,p_8$, 
and a line $L_1$, while any element of $V'$ corresponds to a line $L_2$ through $p_7$ and $p_8$
and a quartic $Q_3$ through the base point of the cubic pencil. For 8 general points we may assume that
the base point of the cubic pencil is not on either conic nor on $L_2$, so
$Q_1+Q_2+L_1=Q_3+L_2$ implies that $L_1$ passes through the base point and 
hence that $L_1\ne L_2$. But this forces $L_2$ to be a component of $Q_1+Q_2$, which it is not.
Thus $V'$ and $W'$ meet only at 0, so $\s(\C G,e_0)=0$, whence $0=\s(\C F_{14},e_0)=\r(\C F_{14},e_0)$.

We now have shown in all cases with $2\le t\le 5$
that $\r(\C F_{6s+t},e_0)\s(\C F_{6s+t},e_0)=0$. 

Turning to $t=1$ with the same method, we find that 
$\s(\C F_{6s+1},e_0)=\s(\C F_{6(s-3)+2},e_0)+\s(\C F_{6s+1}\otimes \C O_E,e_0)$ for all $s\ge 3$,
hence $\s(\C F_{6s+1},e_0)=\s(\C F_{6(s-3)+2},e_0)+16$ for $s\ge 3$, so
$\s(\C F_{6s+1},e_0)=16$ for $s=3$ and for $s\ge 5$, but $\s(\C F_{6s+1},e_0)=17$ for $s=4$.
Using this we check $\r(\C F_{6s+1},e_0)=0$ for $s=3,4$ and $5$, but $\r(\C F_{6s+1},e_0)>0$ for $s>5$.
The case $s=0$ was done above, giving $\s(\C F_1,e_0)=1$ and $\r(\C F_1,e_0)=0$; for $s=1$, consider 
$(0\to\Gamma(\C F_7-\C C)\to\Gamma(\C F_7)\to\Gamma(\C F_7\otimes\C O_C)\to0)\otimes\Gamma(e_0)$
with $\C C=18e_0-7(e_1+e_2+e_3)-6(e_4+\cdots+e_8)$ the union of three $(-1)$-curves, to show
$\r(\C F_7,e_0)=0$ and $\s(\C F_7,e_0)=8$. For $s=2$, consider 
$(0\to\Gamma(\C F_{13}-\C C)\to\Gamma(\C F_{13})\to\Gamma(\C F_{13}\otimes\C O_C)\to0)\otimes\Gamma(e_0)$
with $\C C=36e_0-13(e_1+\cdots+e_6)-12(e_7+e_8)$ the union of six $(-1)$-curves, to show
$\r(\C F_{13},e_0)=0$ and $\s(\C F_{13},e_0)=13$.  Thus we see that 
$\r(\C F_{6s+1},e_0)\s(\C F_{6s+1},e_0)=0$ for $0\le s\le 5$, while $\s(\C F_{6s+1},e_0)=16$ for 
$s>5$.

Finally we have $t=0$; here
$\s(\C F_{6s},e_0)=\s(\C F_{6(s-3)+1},e_0)+\s(\C F_{6s}\otimes \C O_E,e_0)$ for all $s\ge 3$,
hence $\s(\C F_{6s},e_0)=\s(\C F_{6(s-3)+1},e_0)+32$ for $s\ge 3$, from which we get
$\r(\C F_{6s},e_0)=0$ and
$\s(\C F_{6s},e_0)=33, 40, 45, 48, 49, 48$ for $s=3, 4, 5, 6, 7, 8$ respectively,
and $\r(\C F_{6s},e_0)>0$  but $\s(\C F_{6s},e_0)=48$ for $s\ge 9$.
For $s=1$, consider 
$(0\to\Gamma(\C F_6-\C C)\to\Gamma(\C F_6)\to\Gamma(\C F_6\otimes\C O_C)\to0)\otimes\Gamma(e_0)$
with $\C C=18e_0-7(e_1+e_2+e_3)-6(e_4+\cdots+e_8)$ the union of three $(-1)$-curves, to show
$\r(\C F_6,e_0)=0$ and $\s(\C F_6,e_0)=13$. For $s=2$, consider 
$(0\to\Gamma(\C F_{12}-\C C)\to\Gamma(\C F_{12})\to\Gamma(\C F_{12}\otimes\C O_C)\to0)\otimes\Gamma(e_0)$
with $\C C=30e_0-11(e_1+\cdots+e_5)-10(e_6+e_7+e_8)$ the union of five $(-1)$-curves, to show
$\r(\C F_{12},e_0)=0$ and $\s(\C F_{12},e_0)=24$.  Thus we see that 
$\r(\C F_{6s},e_0)\s(\C F_{6s},e_0)=0$ for $0\le s\le 8$, while $\s(\C F_{6s},e_0)=48$ for $s>8$.\qed

We end this section with an example, using our results to obtain a resolution of 
the ideal defining a fat point subscheme.

\rem{Example}{example} Consider eight general points, each
taken with multiplicity $m=205$; thus $Z=205(p_1+\cdots+p_8)$.
Then $I(Z)$ has: 10 generators in degree 579 (since
579 is the first degree $d$ such that $I(Z)_d\ne 0$, and we have
$\hbox{dim}_kI(Z)_{579}=10$); 201 generators in degree
580 (since, for $d=579$, $\C F_d=de_0-205(e_1+\cdots+e_8)$ has 
free part $\C H=51e_0-18(e_1+\cdots+e_8)$
and fixed part $\C N=528e_0-187(e_1+\cdots+e_8)$, and here
$\s(e_0,\C F_d)=\s(e_0,\C H) + (h^0(X,e_0+\C F_d)-h^0(X,\C H+e_0))
=33+168$); 208 in degree 581 (since, for $d=580$, 
$\C H=340e_0-120(e_1+\cdots+e_8)$
and $\C N=240e_0-85(e_1+\cdots+e_8)$, and
$\s(e_0,\C F_d)=48+160$); and 16 in degree 582
(since, for $d=581$, $\C H=581e_0-205(e_1+\cdots+e_8)$
and $\C N=0$, and $\s(e_0,\C F_d)=16+0$). Moreover, the regularity
of $I(Z)$ is 582, so there are no other generators.
(These numbers can be compared with Campanella's bounds
[Cam]: $10\le \nu_{579}\le 10$, $201\le \nu_{580}\le 210$, 
$70\le \nu_{581}\le 280$, and $0\le \nu_{582}\le 79$.)

From this data we easily determine a minimal free resolution for $I(Z)$,
as follows:

$$0\to R^{138}[-581]\oplus R^{216}[-582]\oplus R^{80}[-583]\to
R^{10}[-579]\oplus R^{201}[-580]\oplus R^{208}[-581]\oplus R^{16}[-582]\to I(Z)\to 0.$$

\irSubsection{Nine General Points}{nine}
Now let $p_1,\ldots,p_9$ be distinct
general points of \pr2; clearly, we may assume
that they lie on a smooth plane cubic $C'$. Let 
$X$ be the surface obtained by blowing up the nine
points, and let $C$ denote the proper transform of $C'$.
Note that $C$ is a section of $-K_X$ and that $C$ is numerically effective
with $C^2=0$.

Given any uniform class $\C F=de_0-m(e_1+\cdots+e_9)$, by \ir{recall} we can write
$\C F=(d-3m)e_0-mK_X$, with $h^0(X,\C F)>0$ if and only if $d\ge 3m$, in which case
$h^1(X,\C F)=h^2(X,\C F)=0$, hence $h^0(X,\C F)=(\C F^2-K_X\cdot\C F)/2+1$.
Thus we know $h^0$ for any uniform class \C F. We now determine
$\C S(\C F,e_0)$. 

\prclm{Theorem}{unifnine}{Let $\C F=te_0-mK_X$ with $m\ge 0$,
where $X$ is the blowing up of $r\ge 9$ general points of a smooth plane
cubic $C$.
\item{(a)} If $t>0$ or $t<-1$, then $\C S(\C F,e_0)=0$.
\item{(b)} If $t=0$, then $\s(\C F,e_0)=3m$ and $\r(\C F,e_0)=0$.
\item{(c)} If $t=-1$, then $\s(\C F,e_0)=1$ and $\r(\C F,e_0)=0$.}

\Prf (a) If $t<-1$, then $h^0(X, \C F+e_0)=0$, hence $\C S(\C F,e_0)=0$, so assume that
$t>0$. Note that $-K_X\cdot(te_0-mK_X)>1$
implies $-K_X\cdot(te_0-sK_X)>1$
for all $0\le s\le m$. We will induct on $s$, starting with
the obvious fact that $\C S(te_0,e_0)=0$ for $t\ge 0$. So now we may assume
that $0<s\le m$, and that $\C S(te_0-(s-1)K_X,e_0)=0$. 
Since $h^1(X,te_0-(s-1)K_X)=h^1(X,te_0-(s-1)K_X+e_0)=0$, and
$h^1(X,e_0+K_X)=h^1(X,-e_0)=h^1(\pr2,\C O_{\pr2}(-1))=0$,
by \ir{Mumford}(a,b) we have the exact sequence
$\C S(te_0-(s-1)K_X,e_0)\to\C S(te_0-sK_X,e_0)\to\C S((te_0-sK_X)\otimes\C O_C,e_0\otimes\C O_C)\to 0$,
where the leftmost term vanishes by induction, and $\C S((te_0-sK_X)\otimes\C O_C,e_0\otimes\C O_C)=0$
by \ir{Mumford}(c). Thus $\C S(te_0-sK_X,e_0)=0$ follows by exactness.

(b) Since $h^0(X,-mK_X)=1$, $\C R(-mK_X,e_0)$ clearly vanishes, so
$\s(-mK_X,e_0)=h^0(X,-mK_X+e_0)-h^0(X,-mK_X)h^0(X,e_0)=3m$.

(c) If $t=-1$, then $h^0(X, \C F)=0$ and $h^0(X, \C F+e_0)=1$, so
$\r(\C F,e_0)=0$ and $\s(\C F,e_0)=h^0(X, \C F+e_0)=1$.\qed

\prclm{Corollary}{GIGCnine}{The GIGC holds on \pr2 for $r=9$.}

\Prf By \ir{unifnine}, we see $\r(\C F,e_0)\s(\C F,e_0)=0$ for any uniform class
on $X$. Thus $I(Z)$ has the maximal rank property for any uniform
fat point subscheme $Z$ supported at nine general points of \pr2.\qed

\References

\bibitem{B} Ballico, E. {\it Generators for the homogeneous 
ideal of $s$ general points in \pr 3}, 
J.\ Alg.\  106 (1987), 46--52.

\bibitem{Cam} Campanella, G. {\it Standard bases of perfect homogeneous
polynomial ideals of height $2$}, 
J.\ Alg.\  101 (1986), 47--60.

\bibitem{Cat} Catalisano, M.\ V. {\it ``Fat'' points on a conic}, 
Comm.\ Alg.\  19(8) (1991), 2153--2168.

\bibitem{DGM} Davis, E.\ D., Geramita, A.\ V., and Maroscia, P. {\it Perfect
Homogeneous Ideals: Dubreil's Theorems Revisited},
Bull.\ Sc.\ math., $2^e$ s\'erie, 108 (1984), 143--185.

\bibitem{GO} Geramita, A. V. and Orrechia, F. {\it Minimally
generating ideals defining certain tangent cones}, J.\ Alg. 78
(1982), 36--57.

\bibitem{GM} Geramita, A. V. and Maroscia, P. {\it The ideal
of forms vanishing at a finite set of points of \pr n}, J.\ Alg. 90
(1984), 528--555.

\bibitem{GGR} Geramita, A.\ V., Gregory, D.\ and Roberts, L.
{\it Minimal ideals and points in projective space},
J.\ Pure and Appl.\ Alg. 40 (1986), 33--62.

\bibitem{Gi} Gimigliano, A. {\it Our thin knowledge of fat points},
Queen's papers in Pure and Applied Mathematics, no. 83,
The Curves Seminar at Queen's, vol. VI (1989).

\bibitem{\trans} Harbourne, B. {\it Complete linear systems on rational surfaces}, 
Trans.\ A.\ M.\ S.\ 289 (1985), 213--226. 

\bibitem{\vanc} Harbourne, B. {\it The geometry of rational surfaces and Hilbert
functions of points in the plane},
Can.\ Math.\ Soc.\ Conf.\ Proc.\ 6 
(1986), 95--111.

\bibitem{\ravello} Harbourne, B. {\it Points in Good Position in \pr 2}, in:
Zero-dimensional schemes, Proceedings of the
International Conference held in Ravello, Italy, June 8--13, 1992,
De Gruyter, 1994.

\bibitem{\antican} Harbourne, B. {\it Anticanonical rational surfaces}, preprint
(available via my Web page), 1994.

\bibitem{\fatpts} Harbourne, B. {\it Free Resolutions of Fat Point 
Ideals on \pr2}, preprint (available via my Web page), 1995.

\bibitem{Hi} Hirschowitcz , A.
{\it Une conjecture pour la cohomologie 
des diviseurs sur les surfaces rationelles generiques},
Journ.\ Reine Angew.\ Math. 397
(1989), 208--213.

\bibitem{HSV} Hoa, L.T., Stuckrad, J., and Vogel, W. {\it Towards
a structure theory for projective varieties of degree $=$ codimension $+$ $2$}, 
J. Pure Appl. Algebra 71 (1991), 203--231.

\bibitem{I} Iarrobino, A. {\it Inverse system of a symbolic power, III: thin algebras
and fat points}, preprint 1994.

\bibitem{L} Lorenzini, A. {\it The minimal resolution conjecture}, J.\ Alg.
156 (1991), 5--35.

\bibitem{Ma} Manin, Y.\ I. Cubic Forms. North-Holland Mathematical Library
4, 1986.

\bibitem{Mu} Mumford, D. {\it Varieties defined by quadratric equations},
in: Questions on algebraic varieties, Corso C.I.M.E. 1969 Rome: Cremonese,
1969, 30--100.

\bibitem{N1} Nagata, M. {\it On rational surfaces, II}, Mem.\ Coll.\ Sci.\ 
Univ.\ Kyoto, Ser.\ A Math.\ 33 (1960), 271--293.

\bibitem{N2} Nagata, M. {\it On the 14-th problem of Hilbert}, 
Amer.\ J.\ Math.\ 33 (1959), 766--772.

\bibitem{O} Okuyama, H. {\it A note on conjectures 
of the ideal of $s$-generic points in \pr 4},
J. Math. Tokushima Univ. 25  (1991), 1--11.

\bibitem{Ra} Ramella, I. {\it An algorithmic 
approach to ideal generation of points},
Rend. Accad. Sci. Fis. Mat. Napoli (4) 56 (1989), 71--81.

\bibitem{Ro} Roberts, L. G. {\it A conjecture on Cohen-Macaulay type},
C. R. Math. Rep. Acad. Sci. Canada 3 (1981), no. 1, 43--48.

\bye